\documentclass[amstex,12pt,english,amssymb]{article}

\usepackage{mathtext}
\usepackage[cp1251]{inputenc}
\usepackage[T2A]{fontenc}
\usepackage[english]{babel}
\usepackage[dvips]{graphicx}
\usepackage{amsmath}
\usepackage{amssymb}
\usepackage{amsxtra}
\usepackage{latexsym}
\usepackage{ifthen}

\usepackage[dvips]{graphicx}
\usepackage{epstopdf}
\usepackage{epsfig}
\usepackage{epic}
\usepackage{eepic}

\begin{document}

\newcounter{lemma}[section]
\newcommand{\lemma}{\par \refstepcounter{lemma}%
{\bf Lemma \arabic{section}.\arabic{lemma}.}}
\renewcommand{\thelemma}{\thesection.\arabic{lemma}}

\newcounter{corollary}[section]
\newcommand{\corollary}{\par \refstepcounter{corollary}%
{\bf Corollary \arabic{section}.\arabic{corollary}.}}
\renewcommand{\thecorollary}{\thesection.\arabic{corollary}}

\newcounter{remark}[section]
\newcommand{\remark}{\par \refstepcounter{remark}%
{\bf Remark \arabic{section}.\arabic{remark}.}}
\renewcommand{\theremark}{\thesection.\arabic{remark}}

\newcounter{theorem}[section]
\newcommand{\theorem}{\par \refstepcounter{theorem}%
{\bf Theorem \arabic{section}.\arabic{theorem}.}}
\renewcommand{\thetheorem}{\thesection.\arabic{theorem}}

\newcounter{proposition}[section]
\newcommand{\proposition}{\par \refstepcounter{proposition}%
{\bf Proposition \arabic{section}.\arabic{proposition}.}}
\renewcommand{\theproposition}{\thesection.\arabic{proposition}}

\renewcommand{\theequation}{\arabic{section}.\arabic{equation}}

\def\kohta #1 #2\par{\par\noindent\rlap{#1)}\hskip30pt
\hangindent30pt #2\par}
\def\A{{{\cal {A}}}}
\def\Rn{{{\Bbb R}^n}}
\def\Rk{{{\Bbb R}^k}}
\def\lR{{\overline {{\Bbb R}}}}
\def\lRn{{\overline {{\Bbb R}^n}}}
\def\lRm{{\overline {{\Bbb R}^m}}}
\def\lRk{{\overline {{\Bbb R}^k}}}
\def\lBn{{\overline {{\Bbb B}^n}}}
\def\Bn{{{\Bbb B}^n}}
\let\text=\mbox

\def\Xint#1{\mathchoice
   {\XXint\displaystyle\textstyle{#1}}%
   {\XXint\textstyle\scriptstyle{#1}}%
   {\XXint\scriptstyle\scriptscriptstyle{#1}}%
   {\XXint\scriptscriptstyle\scriptscriptstyle{#1}}%
   \!\int}
\def\XXint#1#2#3{{\setbox0=\hbox{$#1{#2#3}{\int}$}
     \vcenter{\hbox{$#2#3$}}\kern-.5\wd0}}
\def\dashint{\Xint-}

\def\cc{\setcounter{equation}{0}
\setcounter{figure}{0}\setcounter{table}{0}}

\overfullrule=0pt

\title{{\bf On boundary behavior \\ of generalized quasi-isometries}}

\author{{\bf Denis Kovtonyuk  and Vladimir Ryazanov}\\}
\maketitle

\Large \abstract It is established a series of criteria for
continuous and homeomorphic extension to the boundary of the
so-called lower $Q$-homeomorphisms $f$ between domains in
$\overline{\Rn}=\Rn\cup\{\infty\}$, $n\geqslant2$, under integral
constraints of the type $\int\Phi(Q^{n-1}(x))\,dm(x)<\infty$ with a
convex non-decreasing function $\Phi:[0,\infty]\to[0,\infty]$. It is
shown that integral conditions on the function $\Phi$ found by us
are not only sufficient but also necessary for a continuous
extension of $f$ to the boundary. It is given also applications of
the obtained results to the mappings with finite area distortion
and, in particular, to finitely bi-Lipschitz mappings that are a far
reaching generalization of isometries as well as quasi-isometries in
$\Rn$. In particular, it is obtained a generalization and
strengthening of the well-known theorem by Gehring--Martio on a
homeomorphic extension to boundaries of quasiconformal mappings
between QED (quasiextremal distance) domains.
\endabstract
\bigskip



\bigskip
\large \cc
\section{Introduction}
\medskip

In the theory of mappings quasiconformal in the mean, integral
conditions of the type
\begin{equation}\label{eq1.1KR}\int\limits_{D}\Phi(K(x))\,dm(x)<\infty\end{equation}
are applied to various characteristics $K$ of these mappings, see
e.g. \cite{Ah}, \cite{Bi}, \cite{Gol}, \cite{Kr$_1$}--\cite{Ku},
\cite{Per}, \cite{Pes}, \cite{Rya}, \cite{Str}, \cite{UV},
\cite{Zo$_1$} and \cite{Zo$_2$}. Here $dm(x)$ corresponds to the
Lebesgue measure in a domain $D$ in $\Rn$, $n\geqslant2$.
Investigations of classes with the integral conditions
(\ref{eq1.1KR}) are also actual in the connection with the recent
development of the theory of degenerate Beltrami equations, see e.g.
\cite{AIM}, \cite{BGR$_2$}, \cite{BGR$_3$},
\cite{BJ$_1$}--\cite{Dy}, \cite{GMSV}--\cite{IM$_2$}, \cite{Kr$_2$},
\cite{Le}, \cite{MM}--\cite{MRSY}, \cite{MS},
\cite{RSY$_3$}--\cite{SY}, \cite{Tu}, \cite{Ya} and the so-called
mappings with finite distortion, see related references e.g. in the
monographs \cite{IM$_1$} and \cite{MRSY}.

The present paper is a natural continuation of our previous works
\cite{KR$_1$} and \cite{KR$_2$}, see also Chapters 9 and 10 in the
monograph \cite{MRSY}, that have been devoted to integral
conditions of other types turned out to be useful under the study
of mappings with the constraints of the type (\ref{eq1.1KR}).

Recall some definitions. Given a family $\Gamma$ of
$k$-dimensional surfaces $S$ in $\Rn$, $n\geqslant2$,
$k=1,\ldots,n-1$, a Borel function $\varrho:\Rn\to[0,\infty]$ is
called {\bf admissible} for $\Gamma$, abbr.
$\varrho\in\mathrm{adm}\,\Gamma$, if
\begin{equation}\label{eq1.2KR}\int\limits_{S}\varrho^k\,d{\cal{A}}\geqslant1\end{equation}
for every $S\in\Gamma$. The {\bf modulus} of $\Gamma$ is the
quantity \begin{equation}\label{eq1.3KR}M(\Gamma)=
\inf_{\varrho\in\mathrm{adm}\,\Gamma}\int\limits_{\Rn}\varrho^n(x)\,dm(x)\
.\end{equation} We say that a property $P$ holds for {\bf a.e.}
(almost every) $k$-dimensional surface $S$ in a family $\Gamma$ if a
subfamily of all surfaces of $\Gamma$ for which $P$ fails has the
modulus zero.

The following concept was motivated by Gehring's ring definition of
qua\-si\-con\-for\-ma\-li\-ty in \cite{Ge$_3$}. Given domains $D$
and $D'$ in $\lRn=\Rn\cup\{\infty\}$, $n\geqslant2$,
$x_0\in\overline{D}\setminus\{\infty\}$, and a measurable function
$Q:D\to(0,\infty)$, we say that a homeomorphism $f:D\to D'$ is a
{\bf lower Q-homeomorphism at the point} $x_0$ if
\begin{equation}\label{eq1.4KR}M(f\Sigma_{\varepsilon})\geqslant\inf\limits_{\varrho\in\mathrm{adm}\,\Sigma_{\varepsilon}}
\int\limits_{D\cap R_{\varepsilon}}\frac{\varrho^n(x)}{Q(x)}\
dm(x)\end{equation} for every ring
$$R_{\varepsilon}=\{x\in\Rn:\varepsilon<|\,x-x_0|<\varepsilon_0\},\quad\varepsilon\in(0,\varepsilon_0),\
\varepsilon_0\in(0,d_0)\ ,$$ where $$d_0=\sup\limits_{x\in
D}\,|x-x_0|\ ,$$ and $\Sigma_{\varepsilon}$ denotes the family of
all intersections of the spheres
$$S(r)=S(x_0,r)=\{x\in\Rn:|\,x-x_0|=r\},\quad r\in(\varepsilon,\varepsilon_0)\ ,$$
with $D$. The notion can be extended to the case
$x_0=\infty\in\overline{D}$ in the standard way by applying the
inversion $T$ with respect to the unit sphere in $\overline{\Rn}$,
$T(x)=x/|\,x|^2$, $T(\infty)=0$, $T(0)=\infty$. Namely, a
homeomorphism $f:D\to D'$ is a {\bf lower $Q$-homeomorphism at}
$\infty\in\overline{D}$ if $F=f\circ T$ is a lower
Q$_*$-homeomorphism with $Q_*=Q\circ T$ at $0$. We also say that a
homeomorphism $f:D\to{\overline{\Rn}}$ is a {\bf lower
$Q$-homeomorphism in} $D$ if $f$ is a lower $Q$-homeomorphism at
every point $x_{0}\in\overline{D}$.

Further we also give applications of results on lower
$Q$-homeomorphisms to the mappings with finite area distortion
(FAD) and to finitely bi-Lipschitz mappings.

Given domains $D$ and $D'$ in $\Rn$, $n\geqslant2$, following for
\cite{MRSY$_1$}  we say that a ho\-meo\-mor\-phism $f:D\to D'$ is of
{\bf finite metric distortion}, $f\in$ FMD, if $f$ has
$(N)$-property and
\begin{equation}\label{eq1.7KR}0<l(x,f)\leqslant
L(x,f)<\infty\quad{\rm a.e.}\end{equation} where
$$L(x,f)=\limsup_{y\to x}\
\frac{|f(x)-f(y)|}{|x-y|}\ ,$$ and $$l(x,f)=\liminf_{y\to x}\
\frac{|f(x)-f(y)|}{|x-y|}\ .$$

Note that a homeomorphism $f:D\to D'$ is of FMD if and only if $f$
is differentiable with $J(x,f)\neq0$ a.e. and has $(N)$-property,
see Remark 3.11 and Corollary 3.14 in \cite{MRSY$_1$}.

We say that a homeomorphism $f:D\to D'$ has {\bf ($A_k$)-property}
if the two conditions hold:
\medskip

$(A^{(1)}_k)$ : for a.e. $k-$dimensional surface $S$ in $D$ the
restriction $f|_S$ has $(N)$-property with respect to area;
\medskip

$(A^{(2)}_k)$ : for a.e. $k-$dimensional surface $S_*$ in $D'$ the
restriction $f^{-1}|_{S_*}$ has $(N)$-property with respect to
area.
\medskip

We also say that a homeomorphism $f:D\to D'$ is of {\bf finite area
distortion in dimension} $k=1,\ldots,n-1$, $f\in$ FAD$_k$, if $f\in$
FMD and has the $(A_k)-$property. Finally, we say that a
homeomorphism $f:D\to D'$ is of {\bf finite area distortion}, $f\in$
FAD, if $f\in$ FAD$_k$ for every $k=1,\ldots,n-1$. By Lemma 4.1 in
\cite{KR$_1$} every homeomorphism $f\in$ FAD$_{n-1}$ is a lower
$Q$-homeomorphism with $Q(x)$ which is equal to its outer
dilatation. It is known, in particular, that every quasiconformal
mapping belongs to FAD$_{n-1}$, see e.g. Theorem 12.6 in
\cite{MRSY}.

Recall that the {\bf outer dilatation} of a mapping $f:D\to\Rn$,
$n\geqslant2$, at a point $x\in D$ of differentiability for $f$ is
the quantity
$$K_O(x,f)=\frac{||f'(x)||^n}{|J(x,f)|}$$
if $J(x,f)\ne0$, $K_O(x,f)=1$ if $f'(x)=0$, and $K_O(x,f)=\infty$
at the rest points. As usual, here $f'(x)$ denotes the Jacobian
matrix of $f$ at the point $x$, $J(x,f)=\det f'(x)$ is its
determinant and
$$||f'(x)||=\sup\limits_{h\in\Rn\backslash\{0\}}\frac
{|f'(x)h|}{|h|}\ .$$

A homeomorphism $f:D\to D'$ is called {\bf finitely bi-Lipschitz} if
\begin{equation}\label{eq1.8KR}0<l(x,f)\leqslant L(x,f)<\infty\qquad
\forall\ x\in D\ .\end{equation} By Theorem 5.5 in \cite{KR$_1$}
every finitely bi-Lipschitz homeomorphism $f$ is of finite area
distortion and hence it is a lower $Q$-homeomorphism with
$Q(x)=K_O(x,f)$.

\large \cc
\section{Weakly flat and strongly accessible boundaries}
\medskip

Recall first of all the following topological notion. A domain
$D\subset\Rn$, $n\geqslant2$, is said to be {\bf locally
connec\-ted at a point} $x_0\in\partial D$ if, for every
neighborhood $U$ of the point $x_0$, there is a neighborhood
$V\subseteq U$ of $x_0$ such that $V\cap D$ is connected. Note
that every Jordan domain $D$ in $\Rn$ is locally connected at each
point of $\partial D$, see e.g. \cite{Wi}, p. 66.

\begin{figure}[h]
\centerline{\includegraphics[scale=1.0]{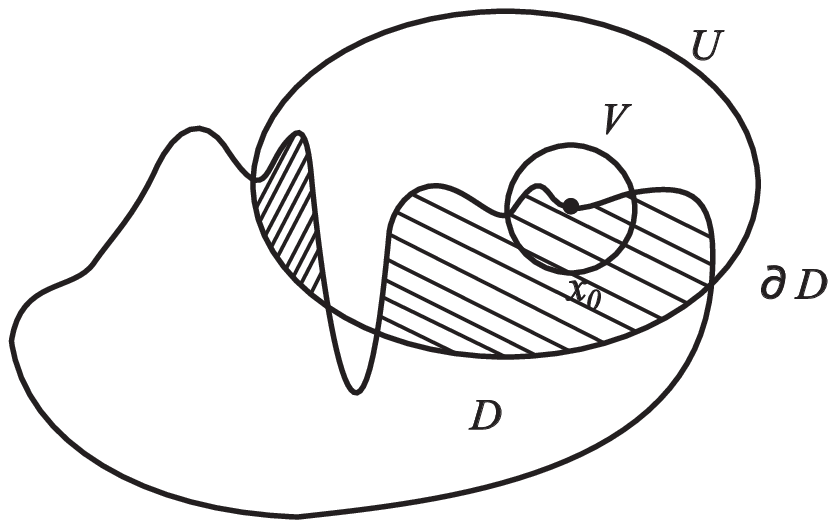}}%
\centerline{Figure 1.}
\end{figure}

We say that $\partial D$ is {\bf weakly flat at a point}
$x_0\in\partial D$ if, for every neighborhood $U$ of the point
$x_0$ and every number $P>0$, there is a neighborhood $V\subset U$
of $x_0$ such that
\begin{equation}\label{eq1.5KR}M(\Delta(E,F;D))\geqslant P\end{equation} for
all continua $E$ and $F$ in $D$ intersecting $\partial U$ and
$\partial V$. Here and later on, $\Delta(E,F;D)$ denotes the
family of all paths $\gamma:[a,b]\to\lRn$ connecting $E$ and $F$
in $D$, i.e. $\gamma(a)\in E$, $\gamma(b)\in F$ and $\gamma(t)\in
D$ for all $t\in(a,b)$. We say that the boundary $\partial D$ is
{\bf weakly flat} if it is weakly flat at every point in $\partial
D$.

\begin{figure}[h]
\centerline{\includegraphics[scale=0.7]{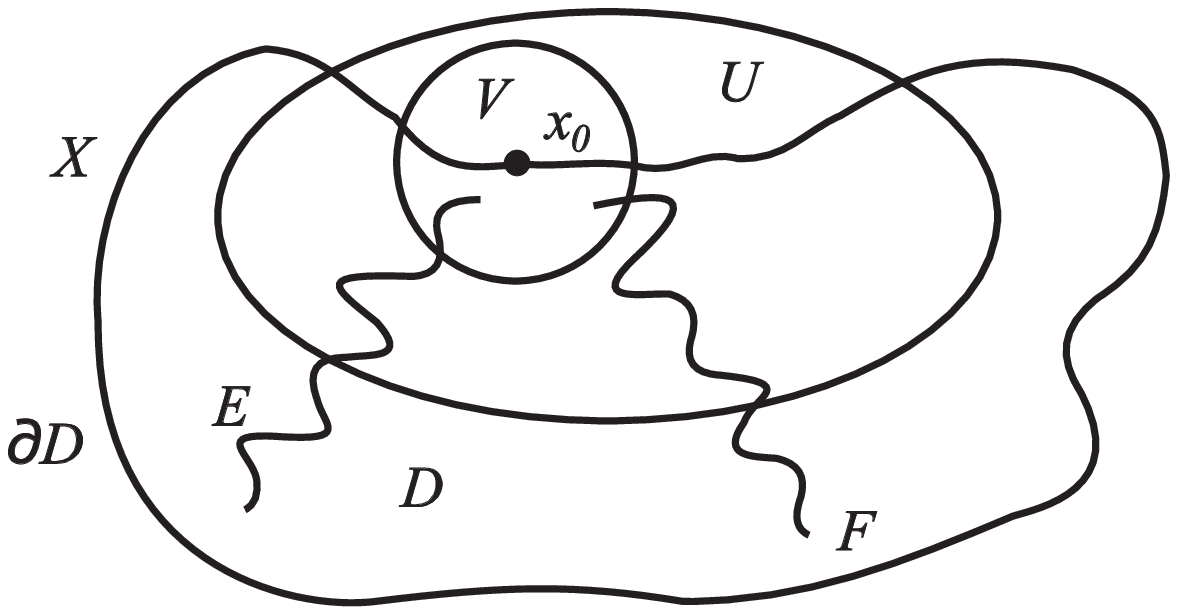}}%
\centerline{Figure 2.}
\end{figure}

We also say that a point $x_0\in\partial D$ is {\bf strongly
accessible} if, for every neighborhood $U$ of the point $x_0$,
there exist a compactum $E$, a neighborhood $V\subset U$ of $x_0$
and a number $\delta>0$ such that
\begin{equation}\label{eq1.6KR}M(\Delta(E,F;D))\geqslant\delta\end{equation} for all
continua $F$ in $D$ intersecting $\partial U$ and $\partial V$. We
say that the boundary $\partial D$ is {\bf strongly accessible} if
every point $x_0\in\partial D$ is strongly accessible.

Here, in the definitions of strongly accessible and weakly flat
boundaries, one can take as neighborhoods $U$ and $V$ of a point
$x_0$ only balls (closed or open) centered at $x_0$ or only
neighborhoods of $x_0$ in another its fundamental system. These
conceptions can also in a natural way be extended to the case of
$\lRn$, $n\geqslant2$, and $x_0=\infty$. Then we must use the
corresponding neighborhoods of $\infty$.

It is easy to see that if a domain $D$ in $\Rn$, $n\geqslant2$, is
weakly flat at a point $x_0\in\partial D$, then the point $x_0$ is
strongly accessible from $D$. Moreover, it was proved by us that
if a domain $D$ in $\Rn$, $n\geqslant2$, is weakly flat at a point
$x_0\in\partial D$, then $D$ is locally connected at $x_0$, see
e.g. Lemma 5.1 in \cite{KR$_2$} or Lemma 3.15 in \cite{MRSY}.

The notions of strong accessibility and weak flatness at boundary
points of a domain in $\Rn$ defined in \cite{KR$_0$} are
localizations and generalizations of the cor\-res\-pon\-ding notions
introduced in \cite{MRSY$_5$}--\cite{MRSY$_6$}, cf. with the
properties $P_1$ and $P_2$ by V\"ais\"al\"a in \cite{Va$_1$} and
also with the quasiconformal accessibility and the quasiconformal
flatness by N\"akki in \cite{Na$_1$}. Many theorems on a
homeomorphic extension to the boundary of quasiconformal mappings
and their generalizations are valid under the condition of weak
flatness of boundaries. The condition of strong accessibility plays
a similar role for a continuous extension of the mappings to the
boundary. In particular, recently we have proved the following
significant statements, see either Theorem 10.1 (Lemma 6.1) in
\cite{KR$_2$} or Theorem 9.8 (Lemma 9.4) in \cite{MRSY}.\bigskip

\begin{proposition}{}\label{prKR2.1} {\it Let $D$ and $D'$ be bounded domains in $\Rn$,
$n\geq2$, $Q:D\to(0,\infty)$ a measurable function and $f:D\to D'$
a lower $Q$-homeomorphism in $D$. Suppose that the domain $D$ is
locally connected on $\partial D$ and that the domain $D'$ has a
(strongly accessible) weakly flat boundary. If
$$\int\limits_{0}^{\delta(x_0)}
\frac{dr}{||\,Q||\,_{n-1}(x_0,r)}=\infty\qquad\forall\
x_0\in\partial D$$ for some $\delta(x_0)\in(0,d(x_0))$ where
$d(x_0)=\sup\limits_{x\in D}\,|x-x_0|$ and
$$||\,Q||\,_{n-1}(x_0,r)=\left(\int\limits_{D\cap
S(x_0,r)}Q^{n-1}(x)\,d\A\right)^\frac{1}{n-1},$$ then $f$ has a
(continuous) homeomorphic extension $\overline{f}$ to $\overline{D}$
that maps $\overline{D}$ (into) onto
$\overline{D'}$.}\end{proposition}\bigskip

Here as usual $S(x_0,r)$ denotes the sphere $|x-x_0|=r$.

\bigskip

A domain $D\subset\Rn$, $n\geqslant2$, is called a {\bf
quasiextremal length domain}, abbr. {\bf QED-domain}, see \cite{GM},
if
\begin{equation}\label{e:7.1}M(\Delta(E,F;\lRn)\leqslant K\cdot
M(\Delta(E,F;D))\end{equation} for some $K\geqslant1$ and all pairs
of nonintersecting continua $E$ and $F$ in $D$.

It is well known, see e.g. \cite{Va$_1$}, that
$$M(\Delta(E,F;\Rn))\geqslant c_n\log{\frac{R}{r}}$$
for any sets $E$ and $F$ in $\Rn$, $n\geqslant2$, intersecting all
the spheres $S(x_0,\rho)$, $\rho\in(r,R)$. Hence a QED-domain has a
weakly flat boundary. One example in \cite{MRSY}, Section 3.8, shows
that the inverse conclusion is not true even among simply connected
plane domains.

A domain $D\subset\Rn$, $n\geqslant2$, is called a {\bf uniform
domain} if each pair of points $x_1$ and $x_2\in D$ can be joined
with a rectifiable curve $\gamma$ in $D$ such that
\begin{equation}\label{e:7.2}s(\gamma)\leqslant a\cdot|x_1-x_2|\end{equation}
and \begin{equation}\label{e:7.3}\min\limits_{i=1,2}\
s(\gamma(x_i,x))\leqslant b\cdot d(x,\partial D)\end{equation} for
all $x\in\gamma$, where $\gamma(x_i,x)$ is the portion of $\gamma$
bounded by $x_i$ and $x$, see~\cite{MaSa}. It is known that every
uniform domain is a QED-domain, but there are QED-domains that are
not uniform, see~\cite{GM}. Bounded convex domains and bounded
domains with smooth boundaries are simple examples of uniform
domains and, consequently, QED-domains as well as domains with
weakly flat boundaries.

A closed set $X\subset\Rn$, $n\geqslant2$, is called a {\bf null-set
of extremal length}, abbr. by {\bf NED-set}, if
\begin{equation}\label{e:8.1} M(\Delta(E,F;\Rn))=M(\Delta(E,F;\Rn\backslash X))\end{equation}
for any two nonintersecting continua $E$ and $F\subset\Rn\backslash
X$.

\begin{remark}\label{rmKR2.0} It is known that if $X\subset\Rn$ is a NED-set, then
\begin{equation}\label{e:8.2}|X|=0\end{equation} and $X$ does not locally separate $\Rn$, i.e.,
\begin{equation}\label{e:8.3}\dim\,X\leqslant n-2\ .\end{equation}
Conversely, if a set $X\subset\Rn$ is closed and
\begin{equation}\label{e:8.4} H^{n-1}(X)=0\ ,\end{equation}
then $X$ is a NED-set, see \cite{Va$_2$}. Note also that the
complement of a NED-set in $\Rn$ is a very particular case of a
QED-domain. \end{remark}

Here $H^{n-1}(X)$ denotes the $(n-1)$-dimensional Hausdorff measure
of a set $X$ in $\Rn$. Also we denote  by $C(X,f)$ the {\bf cluster
set} of the mapping $f:D\to\lRn$ for a set $X\subset\overline D$,
\begin{equation}\label{e:8.5}
C(X,f)\colon=\left\{y\in\overline{{\Bbb
R}^n}:y=\lim\limits_{k\to\infty}f(x_k),\ x_k\to x_0\in X,\ x_k\in
D\right\}.\end{equation} Note, the conclusion $C(\partial
D,f)\subseteq\partial D'$ holds for every homeomorphism $f:D\to D'$,
see e.g. Proposition 13.5 in \cite{MRSY}.

\large \cc
\section{The main lemma}
\medskip

For every non-decreasing function $\Phi:[0,\infty]\to[0,\infty]$,
the {\bf inverse function} $\Phi^{-1}:[0,\infty]\to[0,\infty]$ can
be well defined by setting
\begin{equation}\label{eqKR2.1}\Phi^{-1}(\tau)=\inf\limits_{\Phi(t)\geqslant\tau}t\ .\end{equation}
As usual, here $\inf$ is equal to $\infty$ if the set of
$t\in[0,\infty]$ such that $\Phi(t)\geqslant\tau$ is empty. Note
that the function $\Phi^{-1}$ is non-decreasing, too.

\bigskip

\begin{remark}\label{rmKR2.1} Immediately by the
definition it is evident that \begin{equation}\label{eqKR2.2}
\Phi^{-1}(\Phi(t))\leqslant t\quad\forall\ t\in[0,\infty]
\end{equation} with the equality in (\ref{eqKR2.2}) except
intervals of constancy of the function $\Phi(t)$.\end{remark}

\bigskip

Recall that a function $\Phi:[0,\infty]\to[0,\infty]$ is called
{\bf convex} if $$\Phi(\lambda
t_1+(1-\lambda)t_2)\leqslant\lambda\Phi(t_1)+(1-\lambda)\Phi(t_2)$$
for all $t_1$ and $t_2\in[0,\infty ]$ and $\lambda\in[0,1]$.

In what follows, $\Bn$ denotes the unit ball in the space $\Rn$,
$n\geqslant2$, $$\Bn=\{x\in\Rn:|x|<1\}.$$ The following statement is
a generalization of Lemma 3.1 from \cite{RSY$_1$}.

\bigskip

\begin{lemma}{}\label{lemKR3.1} {\it Let $K:\Bn\to[0,\infty]$ be a measurable
function and let $\Phi:[0,\infty]\to[0,\infty]$ be a
non-decreasing convex function. Then
\begin{equation}\label{eqKR3.1}\int\limits_{0}^{1}\frac{dr}{rk^{\frac{1}{p}}(r)}\geqslant
\frac{1}{n}\int\limits_{eM}^{\infty}\frac{d\tau}{\tau\left[\Phi^{-1}(\tau)\right]^{\frac{1}{p}}}\qquad\qquad\forall\
p\in(0,\infty)\end{equation} where $k(r)$ is the average of the
function $K(x)$ over the sphere $|x|=r$,
\begin{equation}\label{eqKR3.2} M\colon=\dashint_{\Bn}\Phi(K(x))\,dm(x)\end{equation}
is the mean value of the function $\Phi\circ K$ over the unit ball
$\Bn$.}\end{lemma}

\bigskip

\begin{remark}\label{rmKR3.1} Note that (\ref{eqKR3.1}) under every $p\in(0,\infty)$ is
equivalent to \begin{equation}\label{eqKR3.3}
\int\limits_{0}^{1}\frac{dr}{rk^{\frac{1}{p}}(r)}\geqslant
\frac{1}{n}\int\limits_{eM}^{\infty}\frac{d\tau}{\tau\Phi_p^{-1}(\tau)}\qquad\mbox{where}\qquad\
\Phi_p(t)\colon=\Phi(t^p)\ .\end{equation}\end{remark}

{\it Proof of Lemma \ref{lemKR3.1}.} The result is obvious if
$M=\infty$ because then the integral in the right hand side in
(\ref{eqKR3.1}) is zero. Hence we assume further that $M<\infty$.
Moreover, we may also assume that $\Phi(0)>0$ and hence that $M>0$
(the case $\Phi(0)=0$ is reduced to it by approximation of
$\Phi(t)$ through cutting off its graph lower the line
$\tau=\delta>0$). Denote
\begin{equation}\label{eqKR3.4}t_*=\sup\limits_{\Phi_p(t)=\tau_0}
t,\quad\tau_0=\Phi(0)>0.\end{equation} Setting
\begin{equation}\label{eqKR3.5}H_p(t) \colon=\log\,\Phi_p(t),\end{equation} we see that
\begin{equation}\label{eqKR3.6}H_p^{-1}(\eta)=\Phi_p^{-1}(e^{\eta}),\quad\Phi_p^{-1}(\tau)=H_p^{-1}(\log\,\tau).\end{equation}
Thus, we obtain that
\begin{equation}\label{eqKR3.7}k^{\frac{1}{p}}(r)=H_p^{-1}\left(\log\frac{h(r)}{r^n}\right)=
H_p^{-1}\left(n\log\frac{1}{r}+\log\,h(r)\right)\quad\forall\ r\in
R_*\end{equation} where
$h(r)\colon=r^n\Phi(k(r))=r^n\Phi_p\left(k^{\frac{1}{p}}(r)\right)$
and $R_*=\{r\in(0,1):k^{\frac{1}{p}}(r)>t_*\}$. Then also
\begin{equation}\label{eqKR3.8}k^{\frac{1}{p}}(e^{-s})=H_p^{-1}\left(ns+\log\,h(e^{-s})\right)\quad\forall\ s\in S_*\end{equation}
where $S_*=\{s\in(0,\infty):k^{\frac{1}{p}}(e^{-s})>t_*\}$.

Now, by the Jensen inequality and convexity of $\Phi$ we have that
\begin{equation}\label{eqKR3.9}\int\limits_0^{\infty}h(e^{-s})\,ds=\int\limits_{0}^{1}h(r)\,\frac{dr}{r}=
\int\limits_{0}^{1}\Phi(k(r))\,r^{n-1}\,dr\leqslant\end{equation}
$$\leqslant\int\limits_{0}^{1}\left(\dashint_{S(r)}\Phi(K(x))\,d{\cal{A}}\right)\,r^{n-1}{dr}
\leqslant\frac{\Omega_{n}}{\omega_{n-1}}\cdot M=\frac{M}{n}$$
where we use the mean value of the function $\Phi_p\circ K$ over
the sphere $S(r)=\{x\in\Bn:|x|=r\}$ with respect to the area
measure. As usual, here $\Omega_{n}$ and $\omega_{n-1}$ is the
volume of the unit ball and the area of the unit sphere in $\Rn$,
correspondingly. Then arguing by contradiction it is easy to see
that
\begin{equation}\label{eqKR3.10}|T|=\int\limits_{T}ds\leqslant\frac{1}{n}\end{equation}
where $T=\{s\in(0,\infty):h(e^{-s})>M\}$. Next, let us show that
\begin{equation}\label{eqKR3.11}k^{\frac{1}{p}}\left(e^{-s}\right)
\leqslant H_p^{-1}\left(ns+\log\,M\right)\quad\forall\
s\in(0,\infty)\setminus T_*\end{equation} where $T_*=T\cap S_*$.
Note that $(0,\infty)\setminus T_*=\left[(0,\infty)\setminus
S_*\right]\cup\left[(0,\infty)\setminus
T\right]=\left[(0,\infty)\setminus
S_*\right]\cup\left[S_*\setminus T\right]$. The inequality
(\ref{eqKR3.11}) holds for $s\in S_*\setminus T$ by
(\ref{eqKR3.8}) because $H_p^{-1}$ is a non-decreasing function.
Note also that by (\ref{eqKR3.4})
\begin{equation}\label{eqKR3.12}e^{ns}M>\Phi(0)=\tau_0\quad\forall\ s\in(0,\infty)\end{equation}
and then by (\ref{eqKR3.6})
\begin{equation}\label{eqKR3.13}t_*<\Phi_p^{-1}(e^{ns}M)=
H_p^{-1}\left(ns+\log\,M\right)\quad\forall\
s\in(0,\infty).\end{equation} Consequently, (\ref{eqKR3.11}) holds
for $s\in(0,\infty)\setminus S_*$, too. Thus, (\ref{eqKR3.11}) is
true.

Since $H_p^{-1}$ is non-decreasing, we have by (\ref{eqKR3.10})
and (\ref{eqKR3.11}) that \begin{equation}\label{eqKR3.14}
\int\limits_{0}^{1}\frac{dr}{rk^{\frac{1}{p}}(r)}=
\int\limits_{0}^{\infty}\frac{ds}{k^\frac{1}{p}(e^{-s})}\geqslant
\int\limits_{(0,\infty)\setminus
T_*}\frac{ds}{H_p^{-1}(ns+\Delta)}\geqslant\end{equation}
$$\geqslant\int\limits_{|T_*|}^{\infty}\frac{ds}{H_p^{-1}(ns+\Delta)}\geqslant
\int\limits_{\frac{1}{n}}^{\infty}\frac{ds}{H_p^{-1}(ns+\Delta)}=
\frac{1}{n}\int\limits_{1+\Delta}^{\infty}\frac{d\eta}{H_p^{-1}(\eta)}$$
where $\Delta=\log\,M$. Note that $1+\Delta=\log\,eM$. Thus,
\begin{equation}\label{eqKR3.15}\int\limits_{0}^{1}\frac{dr}{rk^{\frac{1}{p}}(r)}\geqslant
\frac{1}{n}\int\limits_{\log\,eM}^{\infty}\frac{d\eta}{H_p^{-1}(\eta)}\end{equation}
and, after the replacement $\eta=\log\,\tau$, we obtain
(\ref{eqKR3.3}), see (\ref{eqKR3.6}), and hence (\ref{eqKR3.1}).

\bigskip

Since the mapping $t\mapsto t^p$ for every positive $p$ is a
sense-preserving ho\-meo\-mor\-phism $[0,\infty]$ onto $[0,\infty]$
we may rewrite Theorem 2.1 from \cite{RSY$_1$} in the following form
which is more convenient for further applications. Here, in
(\ref{eqKR2.4}) and (\ref{eqKR2.5}), we complete the definition of
integrals by $\infty$ if $\Phi_p(t)=\infty$, correspondingly,
$H_p(t)=\infty$, for all $t\geqslant T\in[0,\infty)$. The integral
in (\ref{eqKR2.5}) is understood as the Lebesgue--Stieltjes integral
and the integrals in (\ref{eqKR2.4}) and
(\ref{eqKR2.6})--(\ref{eqKR2.9}) as the ordinary Lebesgue integrals.

\bigskip

\begin{proposition}\label{prKR2.2} {\it Let $\Phi:[0,\infty]\to[0,\infty]$ be a
non-decreasing function. Set
\begin{equation}\label{eqKR2.3}H_p(t)=\log\Phi_p(t),\quad
\Phi_p(t)=\Phi(t^p),\quad p\in(0,\infty).\end{equation} Then the
equality
\begin{equation}\label{eqKR2.4}\int\limits_{\delta}^{\infty}
H^{\prime}_p(t)\,\frac{dt}{t}=\infty\end{equation} implies the
equality \begin{equation}\label{eqKR2.5}
\int\limits_{\delta}^{\infty}\frac{dH_p(t)}{t}=\infty\end{equation}
and (\ref{eqKR2.5}) is equivalent to
\begin{equation}\label{eqKR2.6}\int\limits_{\delta}^{\infty}H_p(t)\,\frac{dt}{t^2}=\infty\end{equation}
for some $\delta>0$, and (\ref{eqKR2.6}) is equivalent to every of
the equalities:
\begin{equation}\label{eqKR2.7}\int\limits_{0}^{\delta}H_p\left(\frac{1}{t}\right)\,dt=\infty\end{equation}
for some $\delta>0$,
\begin{equation}\label{eqKR2.8}\int\limits_{\delta_*}^{\infty}\frac{d\eta}{H_p^{-1}(\eta)}=\infty\end{equation}
for some $\delta_*>H(+0)$,
\begin{equation}\label{eqKR2.9}\int\limits_{\delta_*}^{\infty}\frac{d\tau}{\tau\Phi_p^{-1}(\tau )}=\infty\end{equation}
for some $\delta_*>\Phi(+0)$.

Moreover, (\ref{eqKR2.4}) is equivalent  to (\ref{eqKR2.5}) and
hence (\ref{eqKR2.4})--(\ref{eqKR2.9}) are equivalent each to
other if $\Phi$ is in addition absolutely continuous. In
particular, all the conditions (\ref{eqKR2.4})--(\ref{eqKR2.9})
are equivalent if $\Phi$ is convex and non-decreasing.}
\end{proposition}

\bigskip

It is easy to see that conditions (\ref{eqKR2.4})--(\ref{eqKR2.9})
are more weak under more great $p$, see e.g. (\ref{eqKR2.6}). It
is necessary to give one more explanation. From the right hand
sides in the conditions (\ref{eqKR2.4})--(\ref{eqKR2.9}) we have
in mind $+\infty$. If $\Phi_p(t)=0$ for $t\in[0,t_*]$, then
$H_p(t)=-\infty$ for $t\in[0,t_*]$ and we complete the definition
$H_p'(t)=0$ for $t\in[0,t_*]$. Note, the conditions
(\ref{eqKR2.5}) and (\ref{eqKR2.6}) exclude that $t_*$ belongs to
the interval of integrability because in the contrary case the
left hand sides in (\ref{eqKR2.5}) and (\ref{eqKR2.6}) are either
equal to $-\infty$ or indeterminate. Hence we may assume in
(\ref{eqKR2.4})--(\ref{eqKR2.7}) that $\delta>t_0$,
correspondingly, $\delta<1/t_0$ where
$t_0\colon=\sup\limits_{\Phi_p(t)=0}t$, $t_0=0$ if $\Phi_p(0)>0$.

\large \cc
\section{The main results}
\medskip

Combining Proposition \ref{prKR2.1} and Lemma \ref{lemKR3.1} we
come to the following statement.

\bigskip

\begin{theorem}{}\label{thKR4.1} {\it Let $D$ and $D'$ be bounded domains in $\Rn$, $n\geqslant2$,
$D$ be locally connected on $\partial D$ and $D'$ have (strongly
accessible) weakly flat boundary. Suppose that $f:D\to D'$ is a
lower $Q$-homeomorphism in $D$ with
\begin{equation}\label{eqKR4.1}\int\limits_{D}\Phi(Q^{n-1}(x))\,dm(x)<\infty\end{equation}
for a convex non-decreasing function
$\Phi:[0,\infty]\to[0,\infty]$. If
\begin{equation}\label{eqKR4.2}\int\limits_{\delta_0}^{\infty}\frac{d\tau}{\tau\left[\Phi^{-1}(\tau)\right]^{\frac{1}{n-1}}}=
\infty\end{equation} for some $\delta_0>\tau_0\colon=\Phi(0)$,
then $f$ has a (continuous) homeomorphic extension $\overline{f}$
to $\overline{D}$ that maps $\overline{D}$ (into) onto
$\overline{D'}$.}\end{theorem}

\bigskip

\begin{corollary}{}\label{corKR4.0} {\it If $D$ and $D'$ are either
bounded convex domains or bounded domains with smooth bondaries in
$\Rn$, $n\geqslant2$, then every lower $Q$-homeomorphism $f:D\to D'$
with the conditions (\ref{eqKR4.1}) and (\ref{eqKR4.2}) admits a
homeomorphic extension
$\overline{f}:\overline{D}\to\overline{D'}$.}\end{corollary}

\bigskip

Arguing locally we obtain also the following consequence of
Theorem \ref{thKR4.1}, see Remark \ref{rmKR2.0}.

\bigskip

\begin{corollary}{}\label{corKR4.00} {\it Let $D$ be a domain in $\Rn$,
$n\geqslant2$, and let $X$ be a closed subset of $D$. Suppose that
$f$ is a lower $Q$-homeomorphism of $D\setminus X$ into $\lRn$ such
that
\begin{equation}\label{eqKR4.0}H^{n-1}(X)=0=H^{n-1}(C(X,f))\ .\end{equation}
If the the conditions (\ref{eqKR4.1}) and (\ref{eqKR4.2}) hold,
then $f$ admits a homeomorphic extension to $D$.}\end{corollary}

\bigskip

\begin{remark}\label{rmKR4.1} Note that the condition (\ref{eqKR4.2}) can be rewritten
in the form
\begin{equation}\label{eqKR4.3}\int\limits_{\delta_0}^{\infty}\frac{d\tau}{\tau\Phi_{n-1}^{-1}(\tau)}=
\infty\quad\quad\quad\mbox{where}\quad\Phi_{n-1}(t)\colon=\Phi(t^{n-1})\
.\end{equation} Note also that by Proposition \ref{prKR2.2} the
condition (\ref{eqKR4.3}) can be replaced by every of the condition
(\ref{eqKR2.4})--(\ref{eqKR2.8}) under $p=n-1$ and, in particular,
the condition (\ref{eqKR2.6}) can be rewritten in the form
\begin{equation}\label{eqKR4.4}\int\limits_{\delta}^{\infty}\log\,\Phi(t)\,\frac{dt}{t^{n'}}=+\infty\end{equation}
for some $\delta>0$ where $\frac{1}{n'}+\frac{1}{n}=1$, i.e. $n'=2$
for $n=2$, $n'$ is strictly decreasing in $n$ and $n'=n/(n-1)\to1$
as $n\to\infty$. \end{remark}

\bigskip

{\it Proof of Theorem \ref{rmKR4.1}.} Indeed, let us extend the
function $Q$ by zero outside of $D$ and set, for fixed
$x_0\in\partial D$, $$K(x)=Q^{n-1}(x_0+xd_0)\ ,\qquad x\in\Bn$$ with
some positive $d_0<\sup\limits_{z\in D}|z-x_0|$. Then by Lemma
\ref{lemKR3.1} with the given $K(x)$ and $p=n-1$ we have that
\begin{equation}\label{eqKR4.01}
\int\limits_{0}^{1}\frac{dr}{rk^{\frac{1}{n-1}}(r)}\geqslant
\frac{1}{n}\int\limits_{eM}^{\infty}\frac{d\tau}{\tau\left[\Phi^{-1}(\tau)\right]^{\frac{1}{n-1}}}
\end{equation} where $k(r)$ is the average of $K(x)$ over the
sphere $|x|=r$ and
\begin{equation}\label{eqKR4.02}M=\dashint\limits_{\Bn}\Phi(K(x))\,dm(x).\end{equation}

Now, after the replacement $y_0=x_0+xd_0$ in (\ref{eqKR4.02}), we
have by the condition (\ref{eqKR4.1}) that
$$M\leqslant N\colon=\Phi(0)+\frac{1}{\Omega_nd_0^n}\int\limits_{D}\Phi(Q^{n-1}(y))\,dm(y)<\infty$$
where $\Omega_n$ is the volume of the unit ball in $\Rn$ and after
the replacement $\rho=rd_0$ in the left hand side integral in
(\ref{eqKR4.01}) we obtain that
$$\int\limits_{0}^{d_0}\frac{d\rho}{||Q||_{n-1}(x_0,\rho)}\geqslant
\frac{1}{n\omega_{n-1}^{\frac{1}{n-1}}}\
\int\limits_{eN}^{\infty}\frac{d\tau}{\tau\left[\Phi^{-1}(\tau)\right]^{\frac{1}{n-1}}}$$
where $\omega_{n-1}$ is the area of unit sphere in $\Rn$ and
$$||Q||_{n-1}(x_0,\rho)=\left(\int\limits_{|z-x_0|=\rho}Q^{n-1}(z)\,d{\cal
A}\right)^{\frac{1}{n-1}}.$$

Note that $N>\Phi(0)$. Thus, we conclude from the condition
(\ref{eqKR4.2}) that
\begin{equation}\label{eqKR4.03}\int\limits_{0}^{\delta_0}\frac{d\rho}{||Q||_{n-1}(x_0,\rho)}=\infty.\end{equation}
This is obvious if $\delta\colon=eN\leqslant \delta_0$. If
$\delta>\delta_0$, then
$$\int\limits_{\delta_0}^{\infty}\frac{\delta\tau}{\tau\left[\Phi^{-1}(\tau)\right]^{\frac{1}{n-1}}}=
\int\limits_{\delta}^{\infty}\frac{\delta\tau}{\tau\left[\Phi^{-1}(\tau)\right]^{\frac{1}{n-1}}}+
\int\limits_{\delta_0}^{\delta}\frac{d\tau}{\tau\left[\Phi^{-1}(\tau)\right]^{\frac{1}{n-1}}}$$
where
$$0<\int\limits_{\delta_0}^{\delta}\frac{\delta\tau}{\tau\left[\Phi^{-1}(\tau)\right]^{\frac{1}{n-1}}}\leqslant
\frac{\log{\frac{\delta}{\delta_0}}}{\left[\Phi^{-1}(\delta_0)\right]^{\frac{1}{n-1}}}<\infty$$
because $\Phi^{-1}(\delta_0)>0$.

Finally, by Proposition \ref{prKR2.1} and (\ref{eqKR4.03}) we obtain
the statements of Theorem \ref{rmKR4.1}.

\bigskip

Since quasiconformal mappings are in FAD$_{n-1}$ (of finite area
distortion in dimension $n-1$), see e.g. Theorem 12.6 in
\cite{MRSY}, and QED-domains have weakly flat boundaries, the
following consequence of Theorem \ref{thKR4.1} is a far-reaching
ge\-ne\-ra\-li\-za\-tion of the Gehring-Martio theorem on a
homeomorphic extension to boundaries of quasiconformal mappings
between QED-domains, cf. \cite{GM} and \cite{MV}.

\bigskip

\begin{theorem}{}\label{thKR4.2} {\it Let $D$ and $D'$ be bounded domains in $\Rn$, $n\geqslant2$,
with weakly flat boundaries and let $f:D\to D'$ be a homeomorphism
in FAD$_{n-1}$. If
\begin{equation}\label{eqKR4.5}\int\limits_{D}\Phi(K_{O}^{n-1}(x,f))\,dm(x)<\infty\end{equation}
where $\Phi$ is convex non-decreasing function satisfying at least
one of the conditions (\ref{eqKR2.4})--(\ref{eqKR2.9}) under
$p=n-1$, in particular, (\ref{eqKR4.2}) or (\ref{eqKR4.4}), then $f$
can be extended to a homeomorphism $\overline{f}$ of $\overline{D}$
onto $\overline{D'}$.}\end{theorem}

\bigskip

In turn, since finitely bi-Lipschitz homeomorphisms are of finite
area distortion in dimension $n-1$, see e.g. Theorem 5.5 in
\cite{KR$_1$}, we have also the following consequence.

\bigskip

\begin{corollary}{}\label{corKR4.1} {\it Every finitely be-Lipschitz
homeomorphism $f:D\to D'$ under the hypothesis of Theorem
\ref{thKR4.2} admits a homeomorphic extension to the closure of
the domains $D$ and $D'$.}\end{corollary}

\bigskip

\begin{remark}\label{rmKR4.2} If the domain $D$ is not bounded, then it
must be used the spherical volume $dV(x)=dm(x)/(1+|x|^2)^n$ instead
of the Lebesgue measure $dm(x)$ in the above conditions
(\ref{eqKR4.1}) and (\ref{eqKR4.5}). \end{remark}

\large \cc
\section{Necessary conditions for extension}
\medskip

\begin{theorem}{}\label{thKR5.1} {\it Let $\varphi:[0,\infty]\to[0,\infty]$ be a convex non-decreasing
function such that
\begin{equation}\label{eqKR5.1}\int\limits_{\delta_*}^{\infty}\frac{d\tau}{\tau\varphi^{-1}(\tau)}<\infty\end{equation}
for some $\delta_*\in(\tau_0,\infty)$ where
$\tau_0\colon=\varphi(0)$. Then for every $n\geqslant2$ there is a
diffeomorphism $f$ of the punctured unit ball $\Bn\setminus\{0\}$
onto a ring $\frak{R}=\{x\in\Rn:1<|x|<R\}$ such that
\begin{equation}\label{eqKR5.2}\int\limits_{\Bn}\varphi(K_{O}(x,f))\,dm(x)<\infty\end{equation}
but $f$ cannot be extended by continuity to $0$.}\end{theorem}

\bigskip

By the known criterion of convexity, see e.g. Proposition 5 in I.4.3
of \cite{Bou}, the inclination $[\varphi(t)-\varphi(0)]/t$ is
non-decreasing. By (\ref{eqKR5.1}) the function $\varphi$ cannot be
constant. Thus, the proof of Theorem (\ref{thKR5.1}) is reduced to
the following statement.

\bigskip

\begin{lemma}{}\label{lemKR5.1} {\it Let $\varphi:[0,\infty]\to[0,\infty]$ be a non-decreasing
function such that
\begin{equation}\label{eqKR5.0}\varphi(t)\geqslant C\cdot
t\qquad\forall\ t\in[T,\infty]\end{equation} for some $C>0$ and
$T\in(0,\infty)$ and (\ref{eqKR5.1}) holds. Then for every
$n\geqslant2$ there is a diffeomorphism $f$ of the punctured unit
ball $\Bn\setminus\{0\}$ onto a ring $\frak{R}=\{x\in\Rn:1<|x|<R\}$
such that (\ref{eqKR5.2}) holds but $f$ cannot be extended by
continuity to $0$.}  \end{lemma}

\bigskip

{\it Proof.} Note that by the condition (\ref{eqKR5.1})
\begin{equation}\label{eqKR5.3}\int\limits_{\delta}^{\infty}\frac{d\tau}{\tau\varphi^{-1}(\tau)}<\infty
\quad\quad\quad\forall\ \delta\in(\tau_0,\infty)\end{equation}
because $\varphi^{-1}(\tau)>0$ for all $\tau>\tau_0$ and
$\varphi^{-1}(\tau)$ is non--decreasing. Then applying the linear
transformation $\alpha\varphi+\beta$ with $\alpha=1/C$ and
$\beta=T$, see e.g. (\ref{eqKR2.6}), we may assume that
\begin{equation}\label{eqKR5.4}\varphi(t)\geqslant t\qquad\forall\
t\in[0,\infty)\ .\end{equation} Of course, we may also assume that
$\varphi(t)=t$ for all $t\in[0,1)$ because the values of $\varphi$
in $[0,1)$ give no information on $K_O(x,f)\geqslant1$. It is clear
(\ref{eqKR5.3}) implies that $\varphi(t)<\infty$ for all $t<\infty$,
see the criterion (\ref{eqKR2.6}), cf. (\ref{eqKR2.9}).

Now, note that the function $\Psi(t)\colon=t\varphi(t)$ is strictly
increasing, $\Psi(1)=\varphi(1)$ and $\Psi(t)\to\infty$ as
$t\to\infty$. Hence the functional equation
\begin{equation}\label{eqKR5.5}\Psi(K(r))=\left(\frac{\gamma}{r}\right)^2\quad\quad\quad\forall\
r\in(0,1]\ ,\end{equation} where
$\gamma=\varphi^{1/2}(1)\geqslant1$, is well solvable with $K(1)=1$
and a strictly decreasing continuous $K(r)$, $K(r)<\infty$,
$r\in(0,1]$, and $K(r)\to\infty$ as $r\to0$. Taking the logarithm in
(\ref{eqKR5.5}), we have that
$$\log\,K(r)+\log\,\varphi(K(r))=2\,\log\,\frac{\gamma}{r}$$
and by (\ref{eqKR5.4}) we obtain that
$$\log\,K(r)\ \leqslant\ \log\,\frac{\gamma}{r}\ ,$$
i.e.,
\begin{equation}\label{eqKR5.6}K(r)\ \leqslant\ \frac{\gamma}{r}\ .\end{equation} Then by
(\ref{eqKR5.5}) $$\varphi(K(r))\ \geqslant\ \frac{\gamma}{r}$$ and
by (\ref{eqKR2.2})
\begin{equation}\label{eqKR5.7}K(r)\ \geqslant\ \varphi^{-1}\left(\frac{\gamma}{r}\right)\
.\end{equation}

Next, we define the following mapping in the unit ball $\Bn$:
$$f(x)=\frac{x}{|x|}\,\varrho(|x|)$$
where
$$\varrho(t)\ =\ \exp\{I(t)\},\quad I(t)\ =\ \int\limits_{0}^{t}\frac{dr}{rK(r)}\ .$$
By (\ref{eqKR5.7})
$$I(t)\ =\ \int\limits_{0}^{t}\frac{dr}{rK(r)}\ \leqslant\
\int\limits_{0}^{t}\frac{dr}{r\varphi^{-1}\left(\frac{\gamma}{r}\right)}\
=\
\int\limits_{\frac{\gamma}{t}}^{\infty}\frac{d\tau}{\tau\varphi^{-1}(\tau)}\quad\quad\quad\forall\
t\in(0,1]$$ where $\gamma/t\geqslant\gamma\geqslant1>\varphi(0)=0$.
Hence by the condition (\ref{eqKR5.3})
\begin{equation}\label{eqKR5.8}I(t)\ \leqslant\ I(1)\ =\ \int\limits_{0}^{1}\frac{dr}{rK(r)}\ <\ \infty\quad\quad\quad\forall\
t\in(0,1]\ .\end{equation} Note that $f\in
C^1\left(\Bn\setminus\{0\}\right)$ because $K(r)$ is continuous and,
consequently, $f$ is locally quasiconformal in $\Bn\setminus\{0\}$.

The tangent and radial distortions under the mapping $f$ on the
sphere $|\,x|=\rho$, $\rho\in(0,1)$, are easy calculated
$$\delta_{\tau}(x)\ =\ \frac{|\,f(x)|}{|\,x|}\ =\ \frac{\exp\left\{\int\limits_{0}^{\rho}\frac{dr}{rK(r)}\right\}}
{\varrho}\ ,$$
$$\delta_{r}(x)\ =\ \frac{\partial|\,f(x)|}{\partial|\,x|}\ =\
\frac{\exp\left\{\int\limits_{0}^{\rho}\frac{dr}{rK(r)}\right\}}
{\rho K(\rho)}$$ and we see that
$\delta_{r}(x)\leqslant\delta_{\tau}(x)$ because $K(r)\geqslant1$.
Consequently, by the spherical symmetry we have that
$$K_O(x,f)\ =\
\frac{\delta_{\tau}^{n}(x)}{\delta^{n-1}_{\tau}(x)\cdot\delta_{r}(x)}\
=\ \frac{\delta_{\tau}(x)}{\delta_{r}(x)}\ =\ K(|x|)$$ at all points
$x\in\Bn\setminus\{0\}$, see e.g. Subsection I.4.1 in \cite{Re$_1$}.
Thus, by (\ref{eqKR5.5})
$$\int\limits_{\Bn}\varphi(K_O(x,f))\,dm(x)\ =\ \int\limits_{\Bn}\varphi(K(|x|))\,dm(x)\ =$$
$$=\ \omega_{n-1}\int\limits_{0}^{1}\frac{\Psi(K(r))}{rK(r)}\
r^{n} dr\ \leqslant\
\gamma^{2}\omega_{n-1}\int\limits_{0}^{1}\frac{dr}{rK(r)}\
\leqslant\ M\colon\ =\ \gamma^2\omega_{n-1}I(1)\ <\ \infty\ .$$

On the other hand, along every radial line $x/|x|=\eta\in\Rn
,|\eta|=1$, we have that $f(x)\to\eta$ as $|x|\to0$, i.e. we have no
determinated limit of $f$ under $x\to0$. It is easy to see that
\begin{equation}\label{eqKR5.9}\lim\limits_{x\to0}|f(x)|=\lim\limits_{t\to0}\varrho(t)=e^{0}=1\ ,\end{equation}
i.e. $f$ maps the punctured ball $\Bn\setminus\{0\}$ onto the ring
$1<|y|<R=e^{I(1)}$.

\bigskip

\begin{remark}\label{rmKR5.1} Note that $f$ in the example under the proof of Theorem \ref{thKR5.1}
(Lemma \ref{lemKR5.1}) is finitely bi-Lipschitz and hence of finite
area distortion and, consequently, it is a lower $Q$-homeomorphism
with $Q(x)=K_O(x,f)$, that the domains $D$ and $D'$ have weakly flat
boundaries, the condition (\ref{eqKR5.2}) holds but $f$ cannot be
extended by continuity to the boundary. Thus, the condition
(\ref{eqKR4.2}) in Theorem \ref{thKR4.1} is necessary because if a
function $\Phi$ is convex, then the function $\varphi=\Phi_{n-1}$,
where $\Phi_{n-1}(t)=\Phi(t^{n-1})$, is so. Recall that
$\Phi^{-1}_{n-1}=[\Phi^{-1}]^{\frac{1}{n-1}}$.
\end{remark}

\medskip
\noindent Denis Kovtonyk and Vladimir Ryazanov, \\
Institute of Applied Mathematics and Mechanics,\\
National Academy of Sciences of Ukraine, \\
74 Roze Luxemburg str., 83114 Donetsk, UKRAINE \\
Phone: +38 -- (062) -- 3110145, Fax: +38 -- (062) -- 3110285 \\
denis$\underline{\ \ }$\,kovtonyuk@bk.ru, vl$\underline{\ \
}$\,ryazanov@mail.ru, vlryazanov1@rambler.ru

\end{document}